\documentclass[12pt]{amsart}
\usepackage{amssymb}

\newtheorem{theorem}{Theorem}[section]

\theoremstyle{definition}
\newtheorem{definition}[theorem]{Definition}
\newtheorem{example}[theorem]{Example}

\newtheorem{corollary}[theorem]{Corollary}
\newtheorem{proposition}[theorem]{Proposition}
\theoremstyle{remark}
\newtheorem{remark}[theorem]{Remark}

\begin{document}
\def\C{\mathbb C}
\def\R{\mathbb R}
\def\X{\mathbb X}
\def\cA{\mathcal A}
\def\cT{\mathcal T}
\def\Z{\mathbb Z}
\def\Y{\mathbb Y}
\def\Z{\mathbb Z}
\def\N{\mathbb N}
\def\cal{\mathcal}

\title[Square Root Problem of Kato for the Sum of Operators]{Square Root Problem of Kato for the Sum of Operators}
\author{Toka Diagana} 
\address{Department of Mathematics, Howard University, 2441 6th Street N.W. Wasington D.C. 20059, USA} 
\email{tdiagana@howard.edu}

\begin{abstract}
This paper is concerned with the square root problem of Kato for the
"sum" of linear operators in a Hilbert space ${\mathbb H}$. Under suitable
assumptions, 
we show that if $A$ and $B$ are respectively m-scetroial linear operators satisfying
the square root problem of Kato. Then the same conclusion still holds for
their "sum". As application, we consider perturbed Schr\"{o}dinger operators. 
\end{abstract}
\footnotetext[1]{2000 AMS subject classification. 47A07; 47B44; 47B25.}
\footnotetext[2]{Key words: square root problem of Kato, sesquilinear forms,
m-sectorial operators, algebraic sum, sum form.}
\date{}
\maketitle
\section{Introduction}
In this paper we deal with the square root problem of Kato for the sum of
linear operators in a Hilbert space ${\mathbb H}$. Indeed, let $A, B$ be
(unbounded) m-sectorial operators in a (complex) Hilbert space ${\mathbb H}$ and let
$\Phi$ and $\Psi$ be the (sectorial) sesquilinear forms associated with $A$ and $B$
respectively by the first representation theorem, see, e.g., \cite[Theorem
2.1, p. 322]{kat3}. We say that $A$ and $B$ verify {\it the square root problem
of Kato} if the following holds
\begin{eqnarray}
D(A^{\frac{1}{2}}) = D(\Phi) = D(A^{* \frac{1}{2}}) \; \; \; \mbox{and}
\; \; \; D(B^{\frac{1}{2}}) = D(\Psi) = D(B^{* \frac{1}{2}})
\end{eqnarray}

\medskip

Our primary goal in this paper is to prove that if (1) holds and under suitable
assumptions, then the same conclusion still holds for the algebraic sum
$A + B$, that is,
\begin{eqnarray}
D((A + B)^{\frac{1}{2}}) = D(A^{\frac{1}{2}}) \cap D(B^{\frac{1}{2}}) =
D((A + B)^{* \frac{1}{2}})
\end{eqnarray}
As consequence, we shall discuss the particular case of unbounded normal operators
defined in a (complex) Hilbert space ${\mathbb H}$.

\medskip

It is well-known that the algebraic sum $A+B$ of $A$ and $B$ is not always defined (see \cite{dia2}, \cite{dia3}, and \cite{dia4}). To overcome such
a difficulty, we shall also
deal with an extension of the algebraic sum called {\it sum form}. Recall
that more details about the sum form $A \dotplus B$ of $A$ and
$B$, can be found in \cite[5. Supplementary
remarks, p. 328-32]{kat3} or \cite{dia1}. One then
can show that if (1) holds, and under appropriated assumptions; then the same conclusion still holds for the
sum form, that is,
\begin{eqnarray}
D((A \dotplus B)^{\frac{1}{2}}) = D(A^{\frac{1}{2}}) \cap D(B^{\frac{1}{2}}) = D((A \dotplus B)^{* \frac{1}{2}})
\end{eqnarray}

\medskip

In \cite{mci3},
McIntosh has shown that if $C$ is an invertible m-accretive operator in a
Hilbert space ${\mathbb H}$ such that its spectrum $\sigma(C)$ is a subset
of a region of type $S_{\alpha , \beta} = \{z \in {\mathbb C}: \; \Re e z \geq 0 \;
\mbox{and} \; | \Im m z | \leq \beta (\Re e z )^{\alpha} \}$, where $\alpha \in [0 ,
1)$ and $\beta > 0$.
Then $D(C^{\frac{1}{2}}) = D(C^{* \frac{1}{2}})$. In section 3,
a similar result will be discussed for the sum of invertible m-accretive operators.

\medskip

Historically, the well-known square
problem of Kato takes its origin in a remark formulated in \cite[Remark 2.29,
p. 332-333]{kat3}. It drew the attention of several mathematicians, especially
the pioneer work of McIntosh.

\medskip

Recall that the first counter-example to the square root problem in the
general case of 
abstract m-accretive operators, formulated by Kato,  was found by Lions in \cite{lio2}, that
is, 
\begin{eqnarray}
D(C) = {\mathbb H_{0}^{1}}(0 , + \infty) \; \; \; \mbox{and} \; \; \;
Cu = \displaystyle{\frac{d}{dt}} u, \; \; \; \forall u \in D(C)
\end{eqnarray}
Clearly, $C$ is m-accretive (not m-sectorial) and that: $D(C^{\frac{1}{2}}) \not= D(C^{* \frac{1}{2}})$.

\medskip

A few years later, a remarkable counter-example 
to the square root problem for the general class of abstract m-sectorial operators was found by McIntosh. Indeed, in \cite{mci2},
it is shown that there exists an m-sectorial operator $A$ such that $D(A^{\frac{1}{2}})
\not= D(A^{*\frac{1}{2}})$.
Meanwhile, McIntosh and allies kept investigating on the square root problem
of Kato for elliptic linear operators, formulated by Kato in \cite{kat1}.
Such a question was modified by McIntosh in \cite{mci3}.
Recently, such a famous and challenging question has been solved by McIntosh
and allies. Indeed, they have proven that the domain of the square
root of a uniformly complex elliptic operator 
$A = - div(B\nabla)$ with bounded measurable coefficients in ${\mathbb
R^n}$ is the Sobolev space ${\mathbb H^1}({\mathbb R^n})$ with the 
estimate: $\| A^{\frac{1}{2}} u \|_{L^2} \sim \| \nabla u \|_{L^2}$, where
$\sim$ is the equivalence in the sense of norms, see,
e.g., \cite{aus-tch} and \cite{aus-hof-lac-mci-tch} for details.

\bigskip

\section{Preliminaries}
\subsection{Notation and Definitions}      
Throughout the paper, ${\mathbb R}$, ${\mathbb C}$, $({\mathbb H} , \langle
, \rangle)$, $B({\mathbb H})$ stand for the sets of real, complex numbers, a (complex) Hilbert space endowed with the inner product $\langle , \rangle$
and the space of bounded linear operators, 
respectively; $S_{\alpha , \beta}$ denotes the domain of the complex plan
defined by: $S_{\alpha , \beta} = \{z \in {\mathbb C}: \; \Re e z \geq 0 \;\mbox{and} \; | \Im m z | \leq \beta (\Re e z )^{\alpha} \}$,
where $\alpha \in [0 ,
1)$ and $\beta > 0$. 

\medskip

For a linear operator $A$, 
we denote by $D(A)$, $\sigma(A)$ the domain and the spectrum of $A$. 
For a given sesquilinear form $\Phi: \; D(\phi) \times D(\phi) \subset
{\mathbb H} \times {\mathbb H} \mapsto
{\mathbb C}$, we denote by $\Theta(\phi)$, its numerical range defined
by: $\Theta(\phi) = \{ \phi(u , u): \; u \in D(\phi) \; \; \mbox{with} \; \; \| u \| = 1 \}$. Similarly, the numerical range of a given linear
operator $A$ is defined by:
$\Theta(A) = \{ \langle Au , u \rangle: \; u \in D(A) \; \; \mbox{with} \; \; \| u \| = 1 \}$.

\medskip

Below we list some properties of sectorial sesquilinear forms as well as
m-sectorial operators that we shall use in the sequel.

\medskip

\begin{definition}
A sesquilinear form $\Phi: \; D(\phi) \times D(\phi) \mapsto
{\mathbb C}$ is said to be sectorial if 
$\Theta(\Phi)$ is a subset of the sector of the form
$${\mathcal S_{\alpha , \beta}} = \{ \lambda \in {\mathbb C}: \; 
|\arg( \lambda - \beta) \leq \alpha < {\frac{\pi}{2}}  \},$$
where $\beta \in {\mathbb R}$.
\end{definition}

\medskip

\begin{remark}
Throughout this paper, we assume that $\beta = 0$.  In this case
\begin{eqnarray}
| \Im m \Phi(u , u) | \leq \tan \alpha \; \Re e \Phi (u , u), \; \;
\forall u \in D(\Phi),
\end{eqnarray}
where $\Re e \; \Phi = \displaystyle{\frac{1}{2}}(\Phi + \Phi^*)$
and $\Im m \; \Phi = \displaystyle{\frac{1}{2}}(\Phi - \Phi^*)$ with $\Phi^*$ denotes
the conjugate of the sesquilinear $\Phi$ (see \cite{kat3}).
\end{remark}

\medskip

\begin{definition}
A linear operator $A: \; D(A) \subset {\mathbb H} \mapsto {\mathbb H}$ 
defined on ${\mathbb H}$ is
said to be m-accretive if the following statements hold true
\begin{enumerate}
\item $\Re e \; \langle A u ,u \rangle \; \geq \; 0$
\item $(A + \lambda \; I)^{-1} \in B({\mathbb H})$ and $\|(A +
I \lambda)^{-1} \| \leq \displaystyle{\frac{1}{\Re e \; \lambda}}, \; \; \; \Re e \; \lambda > 0$
\end{enumerate}
\end{definition}

\medskip

\begin{example} Let $\Omega$ be a bounded open subset of ${\mathbb R^n}$
and let $A$ be the operator defined by
$$D(A) = {\mathbb H_{0}^{1}}(\Omega) \cap {\mathbb H^2}(\Omega) \; \; \mbox{with} \; \; A u = - \Delta u,$$
where $\Delta = \sum_{k = 1}^{n} \displaystyle{\frac{\partial^2}{\partial x_{k}^2}}$
denotes the Laplace differential operator.
Clearly, $A$ is (self-adjoint) m-accretive 
in the Hilbert space $L^2(\Omega)$.
\end{example}

\medskip

\begin{definition}A linear operator $A: \; D(A) \subset {\mathbb H} \mapsto {\mathbb
H}$ defined on ${\mathbb H}$ is said to be
quasi-m-accretive if $A+ \xi I$ is m-accretive for some
scalar $\xi$.
\end{definition}

\medskip

\begin{definition}A linear operator $A: \; D(A) \subset {\mathbb H} \mapsto {\mathbb
H}$ defined on ${\mathbb H}$ is said to be
sectorial if $\Theta(A) \subseteq {\mathcal S_{\alpha \; , \; \beta}}$. 
The operator $A$ is said to be m-sectorial if $A$ is sectorial and
quasi-m-accretive.
\end{definition}

\medskip

Let $\Phi$ be a sectorial form in the Hilbert space ${\mathbb H}$.
We denote by ${\mathbb H_{\Phi}}$, the Pre-Hilbert space $D(\Phi)$, when equipped with
the inner product given by 
\begin{eqnarray}
\langle u , v \rangle_{\Phi} \; = \; \Re e \Phi(u , v) + \langle u , v \rangle,
\; \; \; \forall u, v \in D(\Phi)
\end{eqnarray}
It can be shown that ${\mathbb H_{\Phi}}$ is a Hilbert space if and only if
$\Phi$ is a densely defined closed sectorial form.

\medskip

We also need the following theorem due to Lions (see \cite{lio2}).

\medskip

\begin{theorem} Let $A$ be an m-sectorial operator on ${\mathbb H}$ and
let $\Phi$ be the densely defined closed sectorial form associated
with $A$. Assume that there exists a Hilbert space ${\mathbb K} \hookrightarrow {\mathbb
H}$ such that
\begin{enumerate}
\item $D(\Phi)$ is a closed subspace of  $[ {\mathbb K} , {\mathbb H}
]_{\frac{1}{2}}$
\item $D(A) \subset {\mathbb K}$ and  $D(A^*) \subset {\mathbb H}
$
\end{enumerate}
Then 
$$D(A^{\frac{1}{2}})
= D(\Phi) = D(A^{* \frac{1}{2}})$$
\end{theorem}

\medskip

Below we list some properties of the "sum" of operators which we will need in the sequel.

\medskip

\subsection{Sum of Operators}
Let $A, B$ be m-sectorial operators on ${\mathbb H}$.
Their algebraic sum is defined by
$$D(A+B) = D(A) \cap D(B), \; \; \;
(A+B)u  = Au + Bu \; \; \; \forall u \in D(A) \cap D(B)
$$
It is well-known that the algebraic sum defined above is not, always defined.
A typical example can be formulated by the following: Set ${\mathbb H} = L^2({\mathbb R^3})$ and
consider $A, B$, be the m-sectorial
operators given by
$$
D(A) = {\mathbb H^2}({\mathbb R^3}),
\; \; \; Au = - \Delta u, \; \; \forall u \in {\mathbb H^2}({\mathbb R^3})
$$
and
$$
D(B) = \{ u \in L^2({\mathbb R^3}): \; \; V(x) u \in L^2({\mathbb R^3}) \},
\; \; \; B u = Vu, \; \; \forall u \in D(B)
$$
where $V$ is a complex-valued function satisfying the following assumption
\begin{eqnarray}
\Re e \; V > 0, \; \; \; V \in L^1({\mathbb R^3}) \; \; \mbox{and} \; \;
V \not \in L_{loc}^{2}({\mathbb R^3})
\end{eqnarray}

\medskip

\begin{proposition}
Let $A, B$ be the linear operators given above. Assume that the assumption
(7) holds. Then $D(A) \cap D(B) = \{0 \}$.
\end{proposition}

\medskip

\begin{proof}
Let $u \in D(B) \cap D(B)$ and assume that $u \not \equiv 0$.
Since $u \in {\mathbb H^{2}}({\mathbb R^3})$; then $u$ is a continuous function according to the theorem of Sobolev (see \cite{adm}). Thus, there are an open subset $\Omega$ of ${\mathbb R^3}$ and $\delta > 0$ such that $| u(x) | > \delta$ for all $x \in \Omega$.
Let $\Omega '$ be a compact subset of $\Omega$, equipped with the induced topology
by $\Omega$ ($\Omega'$ is also a compact subset of ${\mathbb R^3}$). It follows that,
\begin{eqnarray}
| V |_{\Omega '} = \displaystyle{\frac{ (|Vu|)_{\Omega '}}{|u|_{\Omega '}}} \in L^2(\Omega '),
\end{eqnarray}
Indeed, $(|Vu|)_{\Omega '} \in L^2(\Omega ')$ and $\displaystyle{\frac{1}{(|u|)_{\Omega '}} \in L^{\infty}
(\Omega ')}$. Thus, $V \in L^2(\Omega ')$; this is impossible according to the assumption (7)($V \not \in L_{loc}^{2}({\mathbb R^3})$). Therefore $ u \equiv 0$.
\end{proof}

\medskip

As the previous example shows, the domain of the algebraic sum $A+B$ of 
$A$ and $B$ must be watched carefully. To overcome such a difficulty, we
define an extension of the algebraic sum commonly called sum form, defined with
the help of the sum of sesquilinear forms. Indeed, let $A, B$ be m-sectorial
operators on ${\mathbb H}$ and let $\Phi$ and $\Psi$ be the sesquilinear
forms associated with $A$ and $B$ respectively. It is well-known that
$\Phi$ and $\Psi$ are respectively densely
defined closed sectorial sesquilinear forms. In addition, we have
$$
\Phi (u,v) = \langle Au , v \rangle, \; \; \; \mbox{for every} \; \;
u \in D(A) \; \; \mbox{and} \; \; v \in D(\Phi)\supset D(A)
$$
and
$$
\Psi (u,v) = \langle Bu , v \rangle, \; \; \; \mbox{for every} \; \;
u \in D(B) \; \; \mbox{and} \; \; v \in D(\Psi) \supset D(B)
$$
Now consider their sum defined by, 
$$D(\Xi) = D(\Phi) \cap D(\Psi) \; \; \; \mbox{and} 
\; \; \; \Xi = \Phi + \Psi$$
Assume that 
$\overline{D(\Phi) \cap
D(\Psi)} = {\mathbb H}$; then $\Xi$ is a densely
defined closed sectorial sesquilinear form (see \cite[Theorem 1.31, p.
319]{kat3}). Using the first representation theorem to the 
sectorial sesquilinear
form $\Xi$ (see \cite[Theorem 2.1, p.322]{kat3}); it turns out that there exists a unique m-sectorial operator
associated with it; we denote it by $A \dotplus B$ and call it as
the sum form of $A$ and $B$.

\medskip

Let us notice that the sum form $A \dotplus B$ defined in this way is the
m-accretive extension of the closure $\overline{A + B}$ (if defined) of
$A + B$. Furthermore, $A \dotplus B$ and $\overline{A + B}$ coincide if 
this last is a maximal accretive operator. Therefore, the sum
form $A \dotplus B$ is defined even if $A + B$ is not.

\bigskip

\section{Main Results}
\begin{theorem}Let $A$ and $B$ be m-sectorial
linear operators on ${\mathbb H}$ such that
$$
D(A) = D(A^*) \; \; \; \mbox{and} \; \; \; D(B) = D(B^*)
$$
One supposes that $\overline{D(A) \cap D(B)} =
{\mathbb H}$ and that the closure $\overline{A+B}$ of
$A+B$ is a maximal operator. Then we have
$$
D((\overline{A+B})^{\frac{1}{2}}) = D(A^{\frac{1}{2}}) \cap D(B^{\frac{1}{2}}) = D((\overline{A+B})^{* \frac{1}{2}})
$$
\end{theorem}

\medskip

\begin{proof} Let $\Phi$ and $\Psi$ be the densely defined closed
sectorial sesquilinear forms associated with $A$ and $B$ respectively. Consider
their sum $\Xi = \Phi + \Psi$; since $D(A+B) \subset D(\Phi) \cap
D(\Psi)$ and that $D(A) \cap D(B)$ is dense in $H$. It turns out that
$\Xi$ is a densely defined closed sectorial sesquilinear form
on ${\mathbb H}$. Now, $\overline{A+B}$ is a maximal operator by assumption; it follows that
$\overline{A+B}$ is the operator associated with the sesquilinear form $\Xi$. In the same way, $(\overline{A+B})^*$ is the operator associated with the
conjugate $\Xi^*$  of $\Xi$.
\medskip

Now, $D(A) \cap D(B) = D(A^*) \cap D(B^*)$ with equivalent norms. From the general fact that $A^* + B^* \subset (A+B)^*$.
It follows that $D(\overline{A+B}) \subseteq D((\overline{A+B})^*)$.
Thus,
$$
D((\overline{A+B})^{\frac{1}{2}}) \subseteq D((\overline{A+B})^{* \frac{1}{2}})
$$
Using \cite[Theorem 5.2, p. 238]{lio1}, we obtain that
\begin{eqnarray}
D((\overline{A+B})^{\frac{1}{2}}) \subseteq D(\Xi) \subseteq D((\overline{A+B})^{* \frac{1}{2}})
\end{eqnarray}
Since $\overline{A+B}$ is m-accretive. Then, substituting $\overline{A+B}$ by
$(\overline{A+B})^*$ in (9) yields
\begin{eqnarray}
D((\overline{A+B})^{* \frac{1}{2}}) \subseteq D(\Xi^*) \subseteq D((\overline{A+B})^{\frac{1}{2}})
\end{eqnarray}
Comparing (9) and (10), and using the
fact that $D(\Xi) = D(\Xi^*)$. It follows that,
$
D((\overline{A+B})^{\frac{1}{2}}) = D(A^{\frac{1}{2}}) \cap D(B^{\frac{1}{2}})
= D((\overline{A+B})^{* \frac{1}{2}}).
$
\end{proof}

\medskip

\begin{remark}
Since $A$ and $B$ are respectively m-sectorial; then there $0 \leq \alpha,
\alpha' < \displaystyle{\frac{\pi}{2}}$ such that $\Theta(A) \subset
{\mathcal S_{\alpha , 0}}$ and $\Theta(B) \subset
{\mathcal S_{\alpha' , 0}}$. Setting $\beta = \tan \alpha$ and $\beta' = \tan
\alpha'$; then: 
$$ | \Im m \; \Xi (u,u) | \leq \max(\beta , \beta') \; \Re e \; \Xi(u,u), \; \; \;
\forall u \in D(\Xi) = D(\Phi) \cap D(\Psi)$$
\end{remark}

\medskip

As consequence, we shall apply theorem 3.1 to the case of unbounded normal operators.

\medskip

Let $A$ and $B$ be unbounded normal operators on ${\mathbb H}$.  According
to the spectral theory for unbounded normal operator, one
can write
$$
A = A_{1} - i A_{2} \; \; \mbox{and} \; \; B = B_{1} - i B_{2},
$$
with $A_{k}$, $B_{k}$ self-adjoint operators on ${\mathbb H}$ ($k
= 1, 2$), see, e.g., \cite[pp. 348-355]{rud}. Now since
$D(A) = D(A^*)$ and $D(B) = D(B^*)$, it turns out that
$$
A^* = A_{1} + i A_{2} \; \; \mbox{and} \; \; B^* = B_{1} + i B_{2}
$$
Also, if one supposes that
$A_{k}, B_{k}$ to be nonnegative self-adjoint operators ($k = 1, 2$). Then
$(i A)$ and $(i B)$ are respectively seen as m-accretive operators, see,
e.g., \cite[Corollary 4.4, p. 15]{paz}. Now let us
make the following assumptions

\begin{enumerate}
\item $\; \; \; \overline{D(A) \cap D(B)} = {\mathbb H}$
\item $\exists \; C > 0: \; \;
\langle A_{2} u  , u \rangle\; \leq \; C \; \langle A_{1} u , u \rangle,\; \; \forall
u \in D( A_{1}^{\frac{1}{2}}) \cap D( A_{2}^{\frac{1}{2}})$
\item $\exists \; C' > 0: \; \;
\langle B_{2} u , u \rangle\; \leq \; C' \;  \langle B_{1} u , u \rangle,\; \; \forall
u \in D( B_{1}^{\frac{1}{2}}) \cap D( B_{2}^{\frac{1}{2}})$
\end{enumerate}

\medskip

Here, we set $\Lambda =  D(A^{\frac{1}{2}}) \cap D(B^{\frac{1}{2}})$.

\medskip

\begin{corollary}Let $A = A_{1} - i A_{2}$
and $B = B_{1} - i B_{2}$ be unbounded normal operators on
${\mathbb H}$ such that $A_{k}$ and $B_{k}$ are nonnegative
($k = 1, 2$). Assume that assumptions (i), (ii), and
(iii) hold and that $\overline{A+B}$ is maximal.
Then 
$$
D({\overline{A+B}}^{\frac{1}{2}}) = D(A^{\frac{1}{2}})
\cap D(B^{\frac{1}{2}}) = D({\overline{A+B}}^{* \frac{1}{2}})
$$
\end{corollary}

\medskip

\begin{proof}
Let $\Xi$ the sesquilinear form defined by
$$\Xi(u,v) = \langle (A+B) u , v \rangle,\; \; \forall \; u \in D(A) \cap D(B),
\; \; v \in \Lambda$$
Consider the Pre-Hilbert space ${\mathbb H_{\Xi}} = ( \Lambda\; , \; < , >_{\Xi}$),
where 
$$ \langle u , v \rangle_{\Xi} : =  \langle u , v \rangle_{\mathbb H} + \Re e \Xi (u,v),  \; \; \forall u,
v \in \Lambda$$
Since the sum form operator $A_{1} \dotplus B_{1}$ is a nonnegative
self-adjoint operator. It easily follows that ${\mathbb H_{\Xi}}$ is a Hilbert
space. Thus, the sesquilinear form $\Xi$ is closed. Moreover,
$D(\Xi) = \Lambda$ is dense in ${\mathbb H}$ ($D(A) \cap D(B) \subset \Lambda$
and (i) holds). From the assumptions (ii) and (iii), we conclude that $\Xi$ is
sectorial. Thus, $\Xi$ is a densely defined closed sectorial sesquilinear
form.
According to theorem 3.1, we know that $\overline{A+B}$ is the
m-sectorial operator associated with $\Xi$. Since
$D(A) = D(A^*)$ and $D(B) = D(B^*)$, we complete the proof,
using similar arguments as in the proof of the theorem 3.1.
\end{proof}

\medskip

Let $\Phi$ and $\Psi$ be densely defined closed sectorial
sesquilinear forms on ${\mathbb H}$. Assume that $A$ and $B$ are respectively
the m-sectorial operators associated with $\Phi$ and $\Psi$ by
the first representation theorem (see \cite[Theorem 2.1, p. 322]{kat3}
). Setting $\Xi = \Phi + \Psi$, then we have

\medskip

\begin{theorem} 
Under previous assumptions. One
supposes that $A$, $B$ satisfy (1) and
that $\overline{D(A^{\frac{1}{2}}) \cap D(B^{\frac{1}{2}})} = {\mathbb
H}$. In addition if $D(\Xi)$ is closed in the interpolation space
$[{\mathbb H_{\Xi}} , {\mathbb H}]_{\frac{1}{2}}$. Then there exists
a unique m-sectorial operator $A \dotplus B$ such that
$$D((A \dotplus B)^{\frac{1}{2}}) = D(A^{\frac{1}{2}}) \cap D(B^{\frac{1}{2}})
=  D((A \dotplus B)^{* \frac{1}{2}})$$
\end{theorem}

\medskip

\begin{proof} Since $D(\Xi) = D(A^{\frac{1}{2}}) \cap D(B^{\frac{1}{2}})$
is dense in ${\mathbb H}$. It easily follows that
$\Xi$ is  a densely defined closed sectorial form. According
to Kato's first representation theorem (see \cite[Theorem 2.1. p. 322]{kat3}): there exists a unique
m-sectorial operator $A \dotplus B$ associated with $\Xi$
and that $D(A \dotplus B) \subset D(\Xi) = {\mathbb H_{\Xi}}$,
$D((A \dotplus B)^*) \subset D(\Xi) = {\mathbb H_{\Xi}}$.
Since ${\mathbb H_{\Xi}} \hookrightarrow {\mathbb H}$ is continuous
and that $D(\Xi)$ is closed in $[{\mathbb H_{\Xi}} , {\mathbb H} ]_{\frac{1}{2}}$.
We complete the proof using the theorem of Lions (theorem 2.7).
\end{proof}

\medskip

\begin{theorem}Let $\alpha \in [0 , 1)$ and let $A$ and
$B$ be invertible m-accretive linear operators on
${\mathbb H}$ such that $\overline{D(A) \cap D(B)} = {\mathbb H}$. One supposes that
$\; \Theta(A) \subseteq S_{\alpha, \beta}$ and
$\Theta(B) \subseteq S_{\alpha, \beta'}$, where $\alpha \in [0 , {\frac{\pi}{2}})$
and $\beta, \beta' > 0$. In addition, assume that
$\overline{A+B}$ is m-accretive. Then
\begin{enumerate}       
\item $D((\overline{A+B})^{\frac{1}{2}}) = D((\overline{A+B})^{* \frac{1}{2}}),$
\item $\Theta(\overline{A+B}) \subseteq S_{\alpha, 2\max(\beta, \beta')}$.
\end{enumerate}
\end{theorem}

\medskip

\begin{proof} By assumption $\; \Theta(A) \subseteq S_{\alpha, \beta}$ and
$\Theta(B) \subseteq S_{\alpha, \beta'}$. Thus, we have
\begin{eqnarray}
| \Im m <A u,u> | \leq \beta [ \;  \Re e <A u,u> \; ]^\alpha,\; \;
\forall u \in D(A)
\end{eqnarray}

\begin{eqnarray}
| \Im m <B u,u> | \leq \beta' [ \; \Re e <B u,u> \; ]^\alpha, \; \;
\forall u \in D(B)
\end{eqnarray}
It turns out that, $\forall u \in D(A) \cap D(B)$, there exists $\gamma = \max(\beta,\beta')$ such that,
\begin{eqnarray}
\; \; \; |\Im m < (A+B)u,u > | \leq \gamma [(\Re e <A u,u>)^\alpha +
(\Re e <B u,u>)^\alpha]
\end{eqnarray}
Now note that the following holds: let $\mu \in [0,1]$ and let
$x,y \; \geq 0$. Then
\begin{equation}
x^\mu + y^\mu  \leq 2^{1 - \mu} (x + y )^\mu \leq 2(x + y)^\mu
\end{equation}
Applying "(14)"  to (13), and by density, we have: $\forall u \in
D(\overline{A+B})$

\begin{eqnarray}
| \Im m <\overline{A+B}u,u> | \; \leq 2 \gamma \; [ \Re e < \overline{A+B} u,u> ]^\alpha
\end{eqnarray}

Since $\overline{A+B}$ is m-accretive, we use
(15) and \cite[Theorem B, p. 257-258]{mci3} to
obtain the sought result, that is:
\begin{eqnarray}
D((\overline{A+B})^{\frac{1}{2}}) = D((\overline{A+B})^{*
\frac{1}{2}}
\end{eqnarray}
From (15), it easily follows that $\Theta(\overline{A+B})
\subseteq S_{\alpha, 2\max(\beta,\beta')}.$ 
\end{proof}

\medskip

In what follows, we consider $A, B$ be invertible m-accretive operators
on ${\mathbb H}$ satisfying
\begin{eqnarray}
\Theta(A) \subseteq S_{\alpha , \beta}  \; \; \mbox{and} \; \; \Theta(B) \subseteq S_{\alpha , \beta'},
\end{eqnarray}
where $\alpha \in [0 , 1)$ and $\beta, \beta' > 0$.
Let $\Phi$ and $\Psi$ be the sesquilinear forms
associated with $A$ and $B$, respectively. From (17), it follows that
$A$ and $B$ verify (1). Thus, $\Phi$ and $\Psi$ can be decomposed
as
\begin{eqnarray}
\Phi(u,v) = \langle A^{\frac{1}{2}} u , A^{* \frac{1}{2}} v \rangle,\; \; \; u,v \in D(A^{\frac{1}{2}}) = D(\Phi) = D(A^{* \frac{1}{2}}) ,
\end{eqnarray}
\begin{eqnarray}
\Psi(u,v) = \langle B^{\frac{1}{2}} u , B^{* \frac{1}{2}} v \rangle\; \; \; u,v
\in D(B^{\frac{1}{2}}) = D(\Psi) = D(B^{* \frac{1}{2}}).
\end{eqnarray}
Now consider their sum, $\Xi
= \Phi + \Psi$. Thus, $\forall \; u,v \in D(A^{\frac{1}{2}}) \cap
D(B^{\frac{1}{2}})$,
\begin{eqnarray}
\Xi(u,v) = < A^{\frac{1}{2}} u , A^{* \frac{1}{2}} v > + <
B^{\frac{1}{2}} u , B^{* \frac{1}{2}} v >
\end{eqnarray}
It is not hard to see that $\Theta(\Xi) \subset S_{\alpha , \gamma}$, where
$\alpha \in [0 , 1)$ is given above and $\gamma = 2 \max(\beta , \beta') > 0$. Now, let $A \dotplus B$ be the operator associated with $\Xi$. Thus, we formulate this fact as follows.

\medskip

\begin{theorem} Under previous assumptions; assume $\overline{D(A^{\frac{1}{2}}) \cap D(B^{\frac{1}{2}})} = {\mathbb H}$ and that the operator $A \dotplus
B$ defined above
is (invertible) m-accretive. Then
$$D((A \dotplus B)^{\frac{1}{2}}) = D((A \dotplus B)^{*
\frac{1}{2}} )$$
\end{theorem}

\medskip

\begin{proof}.- Since $A \dotplus B$ is an invertible m-accretive operator
satisfying $\sigma(A \dotplus B) \subset \Theta(A \dotplus B) \subset S_{\alpha , \gamma}$, where
$\alpha \in [0 , 1)$ and $\gamma = 2 \max(\beta , \beta') > 0$.
One completes the proof using a result due to McIntosh \cite[Theorem B, p. 257-258]{mci3}.

\end{proof}

\bigskip

\section{Applications}
This section is concerned with the perturbed Schr\"{o}dinger operators.
Indeed, we shall show that the perturbed operator $S_{Z} = - Z \Delta + V$
verifies the square root problem of Kato, under suitable assumptions on
the complex number $Z$ and the singular complex potential $V$. The
operator $S_{Z}$ will be seen as the algebraic sum of  two m-sectorial
operators $A_{Z}$ and $B$ that we will define in the sequel with
the help of sesquilinear forms.

\medskip

Let $\Omega \subset {\mathbb R^d}$ be an
open subset and set ${\mathbb H} =
L^2(\Omega)$. Let $\Phi_{Z}$ be the sesquilinear form defined by
\begin{eqnarray}
\Phi_{Z}(u,v) = \int_{\Omega} \; Z \nabla u \overline{\nabla v} \; dx, \; \; \;
\; \forall \; u,v \in D(\Phi_{Z}) = {\mathbb H_{0}^{1}}(\Omega),
\end{eqnarray}
where $Z = \alpha - i \beta$ ($\alpha , \beta \in {\mathbb R}$) is a complex number satisfying
\begin{eqnarray}
\alpha, \beta > 0 \; \; \;
\mbox{and} \; \; \; \beta \leq \alpha
\end{eqnarray}
Clearly, the assumption (22) implies that
$\Phi_{Z}$ is a sectorial sesquilinear form on $L^2(\Omega)$.

\medskip

Let $V$ be a measurable complex-valued function and let
$\Psi$ be the sesquilinear form given by
\begin{eqnarray}
\Psi(u,v) = \int_{\Omega} \; V u \overline{v} dx, \; \; \; \; \forall \; u,v
\in D(\Psi),
\end{eqnarray}
where $D(\Psi) = \{u \in L^2(\Omega): \; V | u |^2 \in L^1(\Omega)
\}$. Throughout this section we assume that the potential $V \in
L_{loc}^1(\Omega)$ and that there exists $\theta \in ( 0,
{\frac{\pi}{2}})$ such that
\begin{eqnarray}
| \arg (V(x)) | \leq \theta, \; \; \; \mbox{almost everywhere}
\end{eqnarray}
From (24), it turns out that
\begin{eqnarray}
| \Im m \; \Psi(u,u) | \leq \; \tan \theta \; \Re e \; \Psi(u,u), \; \; \;
\forall \; u \in D(\Psi)
\end{eqnarray}
In other words, $\Psi$ is a sectorial sesquilinear from
on $L^2(\Omega)$. 

\medskip

Under
the previous assumptions, $\Phi$ and $\Psi$ are
respectively, densely defined closed sectorial forms. The operators
associated with both $\Phi_{Z}$
and $\Psi$ are respectively given by

$$D(A_{Z}) = \{u \in {\mathbb H_{0}^{1}}(\Omega) : Z \Delta u \in L^2(\Omega)
\}, \; \; \; A_{Z} u = - Z \Delta u, \; \; \forall \; u \in
D(A_{Z})$$
$$D(B) = \{u \in L^2(\Omega): \; V u \in L^2(\Omega)\} , \; \; B u = V u, \;
\; \forall \; u \in D(B)$$ 
It is not hard to see that $A_{Z}$ and
$B$ are respectively unbounded normal operators on $L^2(\Omega)$ and that
they can be expressed as: $A_{Z} = A_{Z}^{1} - i A_{Z}^2$, where $A_{Z}^1
= - \alpha \Delta$ and $A_{Z}^2 = - \beta \Delta$ are
nonnegative self-adjoint operators, and $B = B_{V}^1 - i
B_{V}^2$, where $B_{V}^1$, $B_{V}^2$ are nonnegative
self-adjoint operators. 

\medskip

Assume that
$\Omega = {\mathbb R^d}$. It will be seen that $\overline{D(A_{Z}) \cap D(B)} = L^2 ({\mathbb R^d})$. Consider the sum $\Xi_{Z} = \Phi_{Z} + \Psi$. Clearly, $\Xi_{Z}$
is a densely defined closed sectorial form. Since $\overline{- Z \Delta + V}$ is m-sectorial (see \cite{bre-kat}). It follows that $\overline{- Z \Delta + V}$ is the operator associated with $\Xi_{Z}$.
In fact, Br\'ezis and Kato computed it in \cite{bre-kat}. It is
defined by
$$D(\overline{- Z \Delta + V}) = \{u \in {\mathbb H^1}({\mathbb R ^d}): \;
V|u|^2 \in L^1({\mathbb R^d}) \; \mbox{and} \; - Z \Delta u + V u \in
L^2({\mathbb R^d}) \}$$
$$\overline{- Z \Delta + V} u = -Z \Delta u + V
u, \; \forall u \in D(\overline{- Z \Delta + V})$$ 
Let us notice
that $D(A_{Z}) = {\mathbb H^2}({\mathbb R^d})$ and $D(B) = \{u \in L^2({\mathbb
R^d}): \; V u \in L^2({\mathbb R^d}) \}$, and their intersection is dense
in $L^2({\mathbb R^d})$. Therefore applying 
Corollary 3.3 to $A_{Z}$ and $B$. It easily follows that
\begin{eqnarray}
D((\overline{- Z \Delta + V})^{\frac{1}{2}}) = {\mathbb H^1}({\mathbb R^d})
\cap D( B^{\frac{1}{2}}) = D((\overline{- Z \Delta + V})^{*
\frac{1}{2}})
\end{eqnarray}
In particular where $d = 1$. Then we obtain that
\begin{eqnarray}
D((\overline{- Z \Delta + V})^{\frac{1}{2}}) = {\mathbb H^1}({\mathbb R})
 = D((\overline{- Z \Delta + V})^{* \frac{1}{2}})
\end{eqnarray}

\bibliographystyle{amsplain}

\end{document}